\documentclass[12pt]{article}
\usepackage{amsfonts}

\title{Generalized eigenvalues for fully tnonlinear singular or degenerate operators in  the radial case.}

\author{F. Demengel}

\date{University of Cergy-Pontoise, 2 avenue Adolphe Chauvain,  95302, Cergy-Pontoise cedex, 
 demengel@math.u-cergy.fr}

\catcode`@=11 \@addtoreset{equation}{section} \catcode`@=12

\newtheorem{theo}{Theorem}[section]
\newtheorem{prop}[theo]{Proposition}
\newtheorem{rema}[theo]{Remark}
\newtheorem{defi}[theo]{Definition}
\newtheorem{cor}[theo]{Corollary}
\newtheorem{lemme}[theo]{Lemma}

\def\R{{\rm I}\!{\rm  R}}
\def\grad{\nabla}

\begin{document} 

\maketitle

\begin{abstract}
In this paper we extend some existence's results concerning the generalized eigenvalues  for  fully nonlinear operators singular or degenerate.  We consider the radial  case and we prove the existence of an infinite number of eigenvalues, simple and isolated. This completes the results  obtained by the author with Isabeau Birindelli  for the first eigenvalues in the  radial case, and the 
results obtained for the Pucci's operator by Busca Esteban and Quaas 
  and for the $p$-Laplace operator  by Del Pino and Manasevich. 
  \end{abstract}
  
 \section{Introduction} 
The extension of the concept of eigenvalue  for   fully nonlinear operators has seen a remarkable  development in these last years, let us mention  the works of  Quaas, Sirakov \cite{QS}, Ishii, Yoshimura \cite{IY}, Juutinen \cite{J}, Patrizi \cite{P},  Armstong \cite{A}, and  previous papers   of the author with Isabeau Birindelli \cite{BD1,BD2} which all deal with the existence of eigenvalues and corresponding eigenfunctions for different fully-nonlinear operators in bounded domains.

\bigskip

In \cite{BD1} we defined the concept of  first eigenvalue on the model of \cite{BNV} and we proved some existence's result  for Dirichlet problem, and for the eigenvalue problem. 

 The  simplicity of the first eigenvalue which is known in the case of the  $p$-Laplacian,  for Pucci's operators, and for operators related but homogeneous  of degree $1$,   remains an open problem for general operators fully nonlinear singular or degenerate  homogeneous  of degree $1+\alpha$ with $\alpha> -1$.  However in \cite{BDr} we proved some uniqueness  result in the case  where  the domain  is  a  ball or an annulus and when  the operator is radial. 
 
 Concerning the "other eigenvalues",   few is known   about  them, except for the Pucci's operators  and  for the $p$-Laplacian,  in the radial case. 
 
 More precisally  in \cite{DM} the authors prove that in the radial case for the $p$-Laplace operator,  there exists an infinite numerable set of eigenvalues, which are simple and isolated,  in \cite{BEQ} the authors  prove the same result for  the Pucci's operators.   Moreover  in each of these papers, the authors   establish some bifurcation results of positive (respectively   negative) solutions for some partial differential equations  related.

 Here we consider also the radial  case for  the model  operator  
 $$F(Du, D^2 u) = |\nabla u|^\alpha{\cal M}_{a,A} (D^2 u)$$
 where  $a$ and $A$ are  two positive numbers, $a\leq A$, $\alpha > -1$  and ${\cal M}_{a,A} $ is   the Pucci's operator 
 ${\cal M}_{a , A} (M) = A tr(M^+)-a tr(M^-)$.
 
  We prove the existence of a numerable set of eigenvalues, $(\mu_k)_k$  which are simple and isolated, and some continuity results for the eigenvalues with respect to the parameters $\alpha, a, A$.

   \section{Assumptions,  notations and previous results in the general case }

   We begin with some generalities about the operators that we consider. 
   
   Let $\Omega$ be some  bounded domain in $\R^N$. 
  
  For $\alpha>-1$ , $F_\alpha$ satisfies :

\begin{itemize}
\item[(H1)] 
 $F_\alpha : \Omega\times \R^N\setminus\{0\}\times S\rightarrow\R$, 
 is continuous and  $\forall t\in \R^\star$, $\mu\geq 0$,
 $F_\alpha (x, tp,\mu X)=|t|^{\alpha}\mu F_\alpha(x, p,X)$.
\item[(H2)]

There exist $0\leq  a\leq A$, such that 
for any $x\in {\Omega}$, $p\in \R^N\backslash \{0\}$, $M\in S$,  $N\in S$, 
$N\geq 0$
\begin{equation}\label{eqaA}
a|p|^\alpha tr(N)\leq F (x,p,M+N)-F (x,p,M) \leq A
|p|^\alpha tr(N).
\end{equation}

\item [(H3)]
 There exists a
continuous function $  \omega$ with $\omega (0) = 0$, such that if
$(X,Y)\in S^2$ and 
$\zeta\in \R^+$ satisfy
$$-\zeta \left(\begin{array}{cc} I&0\\
0&I
\end{array}
\right)\leq \left(\begin{array}{cc}
X&0\\
0&Y
\end{array}\right)\leq 4\zeta \left( \begin{array}{cc}
I&-I\\
-I&I\end{array}\right)$$
and $I$ is the identity matrix in $\R^N$,
then for all  $(x,y)\in \R^N$, $x\neq y$
$$F(x, \zeta(x-y), X)-F(y,  \zeta(x-y), -Y)\leq \omega
(\zeta|x-y|^2).$$
\end{itemize}

Let us now recall the definition of viscosity solutions

\begin{defi}\label{def1}

 Let $\Omega$ be a bounded domain in
$\R^N$, suppose that $f$ is continuous on $\Omega \times \R$,  then
$v$,   continuous in $\Omega$ is called a viscosity super solution (respectively sub-solution)
of
$F(x,\grad u,D^2u)=f(x,u)$ if for all $x_0\in \Omega$, 

-Either there exists an open ball $B(x_0,\delta)$, $\delta>0$  in $\Omega$
on which 
$v= cte= c
$ and 
$0\leq f(x,c)$, for all $x\in B(x_0,\delta)$ (respectively $0\geq f(x,c)$)

-Or
 $\forall \varphi\in {\mathcal C}^2(\Omega)$, such that
$v-\varphi$ has a local minimum on $x_0$ (respectively a local maximum) and $\grad\varphi(x_0)\neq
0$, one has
$$
F( x_0,\grad\varphi(x_0),
 D^2\varphi(x_0))\leq f(x_0,v(x_0)).
$$
(respectively
$$
F( x_0,\grad\varphi(x_0),
 D^2\varphi(x_0))\geq f(x_0,v(x_0)).)
 $$

\end{defi}

 One can also  extend the definition of viscosity solutions to upper semicontinuous sub-solutions and  lower semicontinuous super solutions, as it is done in the paper of Ishii \cite{I}. 
 
  We shall consider in the sequel radial solutions, which will be  solutions of  differential  equations of order two. These solutions will be ${\cal C}^1$ everywhere and ${\cal C}^2$ on each point where their gradient is zero,   so it is easy  to see that these solutions are viscosity solutions.

  We now recall the definition of the  first eigenvalue  and  first eigenfunction adapted to this context, on the model of \cite{BNV}.
  
  We define 
  $$\lambda^+ (\Omega) = \sup \{\lambda, \exists \ \varphi>0,  \  F(x, \nabla \varphi, D^2\varphi)+ \lambda  \varphi^{1+\alpha} \leq 0\ {\rm in} \ \Omega\}$$
    
  $$\lambda^- (\Omega) = \sup \{\lambda, \exists \ \varphi<0,  \  F(x, \nabla \varphi, D^2\varphi)+ \lambda |\varphi|^{\alpha}\varphi  \geq 0\ {\rm in} \ \Omega\}$$
  \begin{rema}
   Let us observe that in this definition,  for $\lambda^+$   (respectively $\lambda^-$), the supremum can be taken over    either continuous and bounded functions, or lower semicontinuous  and bounded  functions (respectively  continuous and bounded functions, or upper semicontinuous and bounded ).
   \end{rema}

We proved in \cite{BD1}  the following existence's result  of "eigenfunctions"
\begin{theo}\label{valp}
 Suppose that $\Omega$ is a bounded regular domain. There exists $\varphi\geq 0$ such that 
 $$\left\{\begin{array} {lc}
 F(x, \nabla \varphi, D^2\varphi)+ \lambda^+(\Omega)  \varphi^{1+\alpha} =0 &\ {\rm in} \ \Omega \\
\varphi=0& {\rm on} \ \partial \Omega
\end{array}\right.$$

 Moreover $\varphi >0$ inside $\Omega$,  is bounded and  continuous. 
  
  Symmetrically  there exists $\varphi\leq 0$ such that 
 $$\left\{\begin{array} {lc}
 F(x, \nabla \varphi, D^2\varphi)+ \lambda^-(\Omega)  |\varphi|^{\alpha}\varphi =0 &\ {\rm in} \ \Omega \\
\varphi=0& {\rm on} \ \partial \Omega
\end{array}\right.$$

 Moreover $\varphi <0$ inside $\Omega$,  is bounded and continuous. 
\end{theo}

These  eigenvalues have the  properties, called maximum and minimum principle   : 

 \begin{theo}
    Suppose that $\Omega$ is a bounded regular domain. 
     If $\lambda < \lambda^+ $,  every  upper semicontinuous and bounded sub-solution of 
     $$F(x,\nabla u, D^2 u) + \lambda |u|^\alpha u \geq 0$$
      which is $\leq 0$ on the boundary, is $\leq 0$ inside $\Omega$.  
  If $\lambda < \lambda^- $,  every  lower semicontinuous and bounded super-solution of 
     $$F(x,\nabla u, D^2 u) + \lambda |u|^\alpha u \leq 0$$
      which is $\geq 0$ on the boundary, is $\geq 0$ inside $\Omega$. 
      \end{theo}
      
      The maximum and minimum principle and some iterative process permit to prove the existence of solutions for the Dirichlet problem,        $$\left\{ \begin{array}{lc}
       F(x,\nabla u, D^2 u) + \lambda |u|^\alpha u = f&{\rm in} \ \Omega\\
        u = 0 & {\rm on} \ \partial \Omega
        \end{array} \right.$$
         where $f$ is supposed to be continuous and bounded, and    $\lambda < \inf (\lambda^+, \lambda^-)$.  
         Moreover if $f\leq 0$ and $\lambda <\lambda^+$ , (respectively $f\geq 0$ and $\lambda <\lambda^-$), there exists a  nonnegative  (respectively non positive) solution .

       We  now give  some  increasing property of  the eigenvalues $\lambda^\pm$  with respect to the domain. 
 
 \begin{prop}\label{propinc}
 
  Suppose that $\Omega$ and $\Omega^\prime$ are some regular  bounded domains such that 
   $\Omega^\prime \subset \subset \Omega$. 
   Then 
   $\lambda^\pm (\Omega^\prime ) > \lambda^\pm (\Omega)$.
    \end{prop}
    
For the convenience of the reader we give a short proof here : 
    We do  it for $\lambda^+$. 
     Let $\varphi$ be an eigenfunction for  $\lambda^+ (\Omega)$. Then by the strict maximum principle there exists $\epsilon >0$ such that 
     $\varphi \geq 2\epsilon$ on $\Omega^\prime $. Define 
     $\lambda^\prime = \lambda^+(\Omega) \inf_{\Omega^\prime}{\varphi^{1+\alpha}\over (\varphi-\epsilon) ^{1+\alpha}}> \lambda^+(\Omega)$. Then the function $\varphi-\epsilon$ is some positive function which satisfies  in $\Omega^\prime$ 
     $$F(x, \nabla (\varphi-\epsilon), \nabla \nabla  (\varphi-\epsilon))+ \lambda^\prime  (\varphi-\epsilon)^{1+\alpha} \leq 0$$
      which implies by the definition of $\lambda^+(\Omega^\prime)$,  that $\lambda^+ (\Omega) < \lambda^\prime \leq \lambda^+(\Omega^\prime)$. 
      
       The following  property of eigenvalues will be needed in section 4   :
       
             \begin{prop}\label{mupropre}
      
       Suppose that there exists $\mu\in \R$, and $u$  continuous and bounded  such that         $$\left\{ \begin{array}{lc}
       F(x,\nabla u, D^2 u)+ \mu |u|^\alpha u=0&,\ u\geq 0, \ u\not \equiv 0 \  {\rm in} \ \Omega\\
       u=0 & {\rm on } \ \partial \Omega
       \end{array}\right. $$
       Then $\mu = \lambda^+$. Symmetrically  
            suppose that there exists $\mu\in \R$, and $u$  continuous  and bounded such that        
            
            $$\left\{ \begin{array}{lc}
             F(x,\nabla u, D^2 u)+ \mu |u|^\alpha u=0&, u\leq 0, \ u\not \equiv 0\ {\rm in} \ \Omega\\
              u= 0 & {\rm on} \ \partial \Omega
              \end{array}\right.$$
              
       Then $\mu = \lambda^-$. 
\end{prop}

 Proof of proposition \ref{mupropre} 
 
 We consider only the first case, the other can be  treated in the same manner. 
 
 By the definition of the first eigenvalue, $\mu\leq \lambda^+$. If $\mu < \lambda^+$, then the minimum principle would imply that $u\leq 0$ in $\Omega$,  a contradiction. 
      
       We  now recall some  regularity and compactness results which will be used in the last section. 
       
        \begin{prop}

Suppose that $\Omega$ is a bounded regular domain. 

Suppose that $F$ satisfies the previous assumptions. Let $f$ be a continuous   and bounded function in ${\Omega}$. Let $u$ be a continuous and bounded viscosity solution of 
\begin{equation}\label{eq4.1}
\left\{
\begin{array}{lc}
F(x, \nabla u, D^2u)=f & \ {\rm in}\
\Omega\\
u=0 &  \ {\rm on}\ \partial\Omega.
\end{array}
\right.
\end{equation}
Then  for any $\gamma<1$ there exists some constant
$C$ which depends only on $|f|_\infty$,  $\gamma$,   $a$, $A$, and $N$,   such that for any
$(x,y)\in\bar\Omega^2$

$$|u(x)-u(y)|\leq C|x-y|^\gamma.$$

\end{prop}

\begin{cor} \label{comp}
Suppose that $\Omega$ is a bounded regular domain. 

Suppose that $F$ satisfies the previous assumptions. Suppose that $(f_n)$ is a  
sequence of continuous and uniformly bounded functions, and $(u_n)$ is a  sequence of  continuous and bounded 
 viscosity solutions of  
$$\left\{\begin{array}{cc}
F( x,  \nabla u_n,\ D^2u_n)=f_n& \ {\rm in} \ \Omega\\
u_n=0 &\ {\rm on} \ \partial \Omega.
\end{array}\right.$$ 
Then the sequence  $(u_n)$
is relatively
compact in
${\mathcal C} (\overline{\Omega})$. 
 Moreover if $f_n $ converges, even simply,  to some continuous and bounded function $f$,  and if for a subsequence $\sigma (n)$,  $u_{\sigma (n)} \rightarrow u$, then $u$ is a solution of the equation with the right hand side $f$. 
\end{cor}

\begin{rema}
Under some additionnal assumption on the regularity of $F$, one has some  Lipschitz regularity  of the solutions.  This assumption is satisfied in the case of the operator considered in the following sections. 
\end{rema} 

We end this section by giving some property of the  first demi-eigenvalues for some particular operators related to Pucci's operators : 

 Let $0<a< A$ and the Pucci's operator 
 $${\cal M}_{a,A} (D^2 u) = A tr( (D^2 u)^+)-a tr((D^2u)^-)$$
  where $(D^2u)^\pm$ denote the positive and negative part of the symmetric matrix $D^2 u$.  
  
   For $\alpha > -1$ the following operator 
   $$F(\nabla u, D^2 u) = |\nabla u|^\alpha {\cal M}_{a,A} (D^2 u)$$
    satisfies the assumption $(H1), (H2)$. We denote by $\lambda_{a,A, \alpha}^\pm$ its  corresponding  first eigenvalues. Then

      \begin{prop} \label{proplam}
       If $a< A$, one has 
       $\lambda_{a,A, \alpha}^+ (\Omega) < \lambda_{a,A, \alpha}^-(\Omega)$.
       
       Moroever if $\lambda_{eq}$ is  the first eigenvalue for the operator $|\nabla u|^\alpha \Delta u$,
        $$\lambda_{a,A, \alpha}^+ \leq a \lambda_{eq} < A \lambda_{eq} \leq \lambda_{a,A, \alpha}^-$$
       \end{prop}
       
       Proof of proposition \ref{proplam}
       
        Let $\phi>0$ be  some  eigenfunction for the eigenvalue $\lambda_{a,A, \alpha}^+ (\Omega)$.
        
         We observe that 
       
                \begin{eqnarray*}
          a \Delta \phi 
           &\leq & A tr(D^2\phi)^+- a tr(D^2\phi)^-\\
           &\leq& {\cal M}_{a, A} (D^2\phi). 
           \end{eqnarray*}
           This implies that 
           $$a \Delta \phi |\nabla \phi|^\alpha + \lambda_{a,A, \alpha}^+ |\phi|^\alpha \phi \leq 0$$
            and then by the definition of $\lambda_{eq}$, 
            $a \lambda_{eq} \geq  \lambda_{a,A, \alpha}^+$. 
            
       In the same manner let $\phi\leq 0$ be such that $\Delta \phi |\nabla \phi|^\alpha = -\lambda_{eq} |\phi|^\alpha\phi$
        then 
        $$|\nabla \phi|^\alpha \left(A tr( (D^2 \phi)^+)- a tr((D^2 \phi)^-) \right)\geq |\nabla \phi|^\alpha  A \Delta \phi = -A \lambda_{eq} |\phi|^\alpha \phi$$
         and by the definition of $\lambda_{a,A, \alpha}^-$ this implies that 
         $$ A \lambda_{eq} \leq      \lambda_{a,A, \alpha}^-.$$
               
  \bigskip
The question of the simplicity of the first eigenvalues for  general operators satisfying $(H1$),.. $(H3)$,   is an open problem. 
The difficulty resides in the fact that one cannot  establish some  strict comparison principle. 
More precisally we should need the following result : 

 {\it If $u\geq v$ and $F(x, \nabla u, D^2 u)= f\leq F(x, \nabla v, D^2 v )= g$ then either $u> v$ everywhere, or $u\equiv v$}.
 
   The difficulty when one wants to prove this result resides on the points where test functions 
   have their gradient equal to  zero.

However we proved in \cite{BDr} the simplicity result   in the radial case.  It  will be precised  in the  forthcoming section, this will be an argument    for the existence and the properties of the other eigenvalues in the case of the operator $|\nabla u|^\alpha {\cal M}_{a,A} (D^2 u)$.

 \section{The radial case}

   Let $\Omega$ be a ball $B(0,1)$ or an annulus  $B(0,1)\setminus \overline{B(0, \rho)}$ for some $\rho\in ]0, 1[$. 
 
 We suppose   that there exists $\tilde F$  such that for any radial function $u(x) = g(|x|)$, 
$F(x, \nabla u, D^2 u) = \tilde F(r, g^\prime , g^{\prime\prime})$. 
 In that case the conditions on $F$ imply that 
 $$ |g^\prime|^\alpha\left(\gamma_1 g^{\prime\prime}+\frac{\gamma_2(N-1)}{|x|}g^\prime\right)\leq F(x,\grad \phi,D^2\phi)\leq  |g^\prime|^\alpha\left(\Gamma_1 g^{\prime\prime}+\frac{\Gamma_2(N-1)}{|x|}g^\prime\right)$$
where 
$$\gamma_1=\left\{\begin{array}{lc}
                   a & {\rm if}\  g^{\prime\prime}>0\\
                   A &   {\rm if}\ g^{\prime\prime}<0
                   \end{array}, \right. \  \gamma_2=\left\{\begin{array}{lc}
                   a &  {\rm if}\  g^{\prime}>0\\
                   A &  {\rm if}\  g^{\prime}<0, \end{array} \right.\ 
                   $$
                   $$
 \Gamma_1=\left\{\begin{array}{lc}
                   A &   {\rm if}\  g^{\prime\prime}>0\\
                   a &  {\rm if}\  g^{\prime\prime}<0
                   \end{array}, \right. \  \Gamma_2=\left\{\begin{array}{lc}
                   A &  {\rm if}\  g^{\prime}>0\\
                   a &  {\rm if}\  g^{\prime}<0. \end{array} \right.
                   $$

     In this situation one can define the first  radial eigenvalues $\lambda^\pm_{rad} (\Omega)$
    
    $$\lambda_{rad}^+ (\Omega ) = \sup \{\lambda , \exists \varphi>0, \  {\rm radial},\ 
     \tilde F(r, \varphi^\prime , \varphi^") + \lambda |\varphi|^\alpha \varphi\leq 0\ {\rm in} \ \Omega \}$$
      $$\lambda_{rad}^- (\Omega ) = \sup \{\lambda , \exists \varphi<0, \ {\rm radial},  \ 
     \tilde F(r, \varphi^\prime , \varphi^") + \lambda |\varphi|^\alpha \varphi\geq 0\ {\rm in} \ \Omega \}$$
     
      Acting as in the general case, one can prove the existence of eigenfunctions for each of these eigenvalues, and using the maximum and minimum principle one  derives that  $\lambda^\pm _{rad} (\Omega )= \lambda^\pm (\Omega)$ in the sense given in theorem   \ref{valp} for the operator $F(x, \nabla u, D^2 u)$.

\begin{rema}\label{remaru}
  
  In the case of  the ball,   for any  constant sign viscosity   solution   of 
 $$\left\{\begin{array}{lc}
 \tilde F(r, u^\prime, u^{\prime\prime})  + \lambda^\pm  |u|^\alpha u = 0\ &{\rm in }  \ B(0,1)\\
  u = 0 &{\rm on} \ \{r = 1\},
  \end{array}\right. $$
  Then $u$ is decreasing from $r=0$  for $\lambda^+$,  increasing from $r=0$ for $\lambda^-$. 
  In particular if $u$ is ${\cal C}^1$,  $0$ is the unique point where $u^\prime$ is zero.

  In the case of an annulus $B(0,1)\setminus \overline{B(0, \rho)}$, if $u$ is a positive (respectively negative)   viscosity  solution 
 of 
   \begin{equation}\label{equaz}
   \left\{\begin{array}{lc}
 \tilde F(r, u^\prime, u^{\prime\prime})   + \lambda^\pm  |u|^\alpha u = 0\ &{\rm in }  \  B(0,1)\setminus \overline{B(0, \rho)}\\
  u = 0 &{\rm on} \ \{r = 1\}\ {\rm and} \ \{r = \rho\},
  \end{array}\right. \end{equation}
  then  there exists a unique point $r = r_u$ such that  
  $u$ is increasing (respectively decreasing) on $[\rho, r_u]$,  and decreasing (respectively increasing ) on $[r_u, 1]$. 
  In particular if $u$ is ${\cal C}^1$,  $r_u$ is the unique point where $u^\prime$ is zero.  
  
  \end{rema}
  The uniqueness result obtained in \cite{BDr} is the following : 
  
         \begin{prop}\label{propun}
        Suppose that $\Omega$ is a ball or an annulus. 
        Suppose that $\varphi$ and $\psi$ are two  positive  radial eigenfunctions  in the viscosity sense,  for the eigenvalue $\lambda^+$,  which are zero on the boundary,  then there exists some  positive constant $c$ such that 
        $\varphi = c\psi$.
        \end{prop}
        
     \begin{rema}
     Of  course the same result holds for  the negative eigenfunctions  corresponding to  $\lambda^-$.
     \end{rema}
   \bigskip

    From now we shall  denote by an abuse of notation by 
    ${\cal M}_{a,A} (r, g^\prime, g")$ the operator $g\mapsto \Gamma_1 g^{\prime\prime}+\frac{\Gamma_2(N-1)}{r}g^\prime$  and $\tilde F$ will be 
    \begin{equation}\label{eqpucci}
    \tilde F(r,  g^\prime , g^") =   |g^\prime|^\alpha\left(\Gamma_1 g^{\prime\prime}+\frac{\Gamma_2(N-1)}{r}g^\prime\right), 
    \end{equation}
    where $\Gamma_1$ and $\Gamma_2$ are the multivalued functions defined at the beginning of section 3. 
    \begin{rema}
    We shall most of the time use  more correctly   the  definition which is valid when $g$ is Lipshitz,   and when $\Gamma_1$ and $\Gamma_2$ are  determined :       $$\tilde F (r, g^\prime, g^")= \Gamma_1 {d\over dr} ({|g^\prime|^\alpha g^\prime \over 1+\alpha}) + \Gamma_2 {(N-1)\over r} |g^\prime |^\alpha g^\prime,  $$
    the derivative ${d\over dr} ({|g^\prime|^\alpha g^\prime) \over 1+\alpha})$ being taken in the distributional sense.  
     
     \end{rema} 
     
     We end this section by giving  one consequence of the Hopf principle in the case of  the operator $\tilde F$. 
          
                           \begin{rema}\label{remhopf}
                           Suppose that  $u$  is a non negative  solution  in the viscosity sense of $\tilde F(r, u^\prime, u^") = f$   on $[0, R[$ for some $R \leq \infty$, with $f$ continuous and non positive ,   then either $u>0$ everywhere, or $u\equiv0$; In particular if $u$ satisfies      $\tilde F(r, u^\prime, u^") = -\lambda |u|^\alpha u$ with $\lambda >0$,  and if  $u(r_o)=0$ then $u$ must change sign on $r_o$.                         
                            \end{rema}

    \section{ The functions $w^+$ and $w^-$. }
    
    In this section we prove the existence  and uniqueness of some radial  solutions of 
  $$ \left\{ \begin{array}{lc}
   |w^\prime |^\alpha {\cal M}_{a, A} (r, w^\prime, w^") = -|w|^\alpha w&\  {\rm in} \ \R^+\ , \\
   w(0) = 1,  w^\prime (0) = 0&\ 
   \end{array}\right.$$
   This will permit as in \cite {BEQ}, \cite{DM} to prove  the existence of an infinite numerable set of  radial eigenvalues for the operator $|\nabla w |^\alpha {\cal M}_{a,A} (D^2 w)$ in the ball.

\begin{prop}\label{propex}
 There exists a unique ${\cal C}^1$ solution  of the equation 
 \begin{equation} \label{eq1}
 |w^\prime |^\alpha ({\cal M}_{a, A} (r, w^\prime, w^") )= -|w|^\alpha w\  {\rm in} \ \R^+\ , 
  w(0) = 1,  w^\prime (0) = 0
  \end{equation}
  Moreover $w$ is ${\cal C}^2$  around  each point where $ w^\prime \neq 0$. 
 
\end{prop}

This proposition will be a consequence of the three following results :

 \begin{prop}\label{propk}
 For all $r_o\geq 0$,  and for all $k_o\neq 0$
 there exists some $\delta >0$ such that there  is existence and uniqueness of solution to 
  \begin{eqnarray} \label{eq2}
&& a \left( |k^\prime |^\alpha k^\prime\left( {N-1\over r}\right) + {d\over dr} ({ |k^\prime |^\alpha k^\prime\over 1+\alpha} )\right) = -|k|^\alpha k\  {\rm for } \  r\in ] r_o, r_o+\delta[\ ,  \ {\rm or}\  r\in ]r_o-\delta, r_o[\cap \R^+, \nonumber\\
 && k(r_o) = k_o,  k^\prime (r_o) = 0,\       
  \end{eqnarray}
 \begin{eqnarray} \label{eq3}
  && A|k^\prime |^\alpha k^\prime\left( {N-1\over r}\right) +a {d\over dr} ({ |k^\prime |^\alpha k^\prime\over 1+\alpha} ) = -|k|^\alpha k\  {\rm for } \  {\rm for} \  r\in ] r_o, r_o+\delta[\ , {\rm or}  \  r\in ] r_o-\delta, r_o[\cap \R^+\ ,
   \nonumber\\
  &&k(r_o) = k_o,  k   ^\prime (r_o) = 0,
  \end{eqnarray}
    \begin{eqnarray} \label{eq4}
&& A \left(|k^\prime |^\alpha k^\prime \left({N-1\over r} \right)+{d\over dr} ({ |k^\prime |^\alpha k^\prime\over 1+\alpha} )\right) = -|k|^\alpha k\  {\rm for } \  r\in ] r_o, r_o+\delta[\ ,  \ {\rm or}\  r\in ]r_o-\delta, r_o[\cap \R^+,\nonumber\\
 &&   k(r_o)=  k_o,  k^\prime (r_o) =0,
  \end{eqnarray}
      \begin{eqnarray} \label{eq5}
&& a|k^\prime |^\alpha k^\prime \left({N-1\over r} \right)+A {d\over dr} ({ |k^\prime |^\alpha k^\prime\over 1+\alpha} ) = -|k|^\alpha k\  {\rm for } \ \   r\in ] r_o, r_o+\delta[\ ,  \ {\rm or}  \  r\in ] r_o-\delta, r_o[\cap \R^+  \nonumber\\ 
&&  k(r_o) =k_o,  k^\prime (r_o) = 0.
  \end{eqnarray}
   Moreover $k$ is ${\cal C}^2$  around  each point where $ k^\prime \neq 0$. 
 
  \end{prop}
  
    In a second step we shall prove the existence's and uniqueness result : 
    
     \begin{prop}\label{propmM}
   If $w_o^\prime\neq 0$,  and for all $w_o$, there exists a local  unique solution  to 
   $$  {\cal M}_{a,A} (r, w^\prime,  w^") = -{|w|^\alpha w\over |w^\prime |^\alpha}$$ 
   $$ (w(r_o), w^\prime (r_o)) = (w_o, w_o^\prime) $$
   Moreover   if  on $]r_1, r_2[\subset ]0,\infty[$,  $w$ is a maximal  solution,     $\lim_{r\rightarrow  r_i, \ r\in ]r_1, r_2[} w^\prime (r) = 0$,  $w$ is ${\cal C}^2$ on $]r_1, r_2[$, $ {d\over dr} (|w^\prime |^\alpha w^\prime(r)) $  exists everywhere  on $]r_1, r_2[$ and 
   ${d\over dr} (|w^\prime |^\alpha w^\prime)(r_1^+)w(r_1^+)<0$,  and  ${d\over dr} (|w^\prime |^\alpha w^\prime)(r_2^-)w(r_2^-)<0$. 
    
    \end{prop}
        
      \bigskip
   \begin{prop}\label{propr1}
    Let $\delta$ be such that on ${\cal C}([0, \delta])$,  $k$ in (\ref{eq2})  with $r_o = 0$ and $k_o = 1$  is well defined and $|k-1|_{{\cal C} ([0, \delta])}  < {1\over 2} $. 
    Then   there exists some constant $c_1$  depending on $a$, $A$, $N$ such that 
    $|k^\prime |\leq c_1$. Moreover there exists $r_1 >0$ which depends only on $a$, $A$, and $N$ such that $k^\prime$ and $k^"$ are $<0$ on $]0, r_1[$. 
    
    \end{prop}
    
     \begin{rema}
     The analogous result holds for the  situations in (\ref{eq3}), (\ref{eq4}), (\ref{eq5}). 
     \end{rema}
           
      We postpone the proof of  these three propositions, and we    conclude to the local  existence and uniqueness's  result, arguing as follows  : 
   
  Let $r_o = 0$, $k$ be the solution of (\ref{eq2})  with $k_o = 1$, and, according to proposition \ref{propr1},  let $r_1$ be such that on $]0, r_1]$,   $k^\prime$ and $k^{\prime\prime}$ are negative. Let $w$ be the solution  given by proposition \ref{propmM} of 
      \begin{equation} \label{eq6}
{\cal M}_{a, A} (r, w^\prime, w^{\prime\prime}) = -{|w|^\alpha w\over  |w^\prime |^\alpha }\  {\rm in} \ \R^+\ , 
  w(r_1) = k(r_1), w^\prime (r_1) =k^\prime (r_1)\neq 0
  \end{equation}on some neighborhood 
$]r_1-\delta_1, r_1[$.   By the equation one must have $w^{\prime\prime}(r_1) <0$. Then by uniqueness $  w= k$ on $]r_1-\delta_1, r_1[$.  We can continue replacing $r_1$ by $r_1-\delta_1$ and finally obtain that $w = k$ on the left of $r_1$ as long as $w^\prime\neq 0$, i.e. until $0$.  So we have obtained the existence and uniqueness of solution on a neighborhood  on the right of zero. 
  
   We can extend the solution on the right of $r_1$. If $w^\prime (r) \neq 0$ for all $r \geq r_1$,  the result is given by proposition \ref{propmM}.  Suppose  now that $r_o\geq  r_1$ is  the first point after $r_1$ such  that $w^\prime (r_o) = 0$.  By remark \ref{remhopf} in section 3, $w(r_o)$ cannot be zero. If $w(r_o) <0$,   anticipating on the behaviour of the possible solutions on the right of $r_o$, we know by  using  the conclusion in proposition \ref{propmM},  that one must have $\lim_{r\rightarrow r_o, r> r_o} {d\over dr} (|w^\prime |^\alpha w^\prime(r)) >0$ , so the equation to solve on the right of $r_o$  is   (\ref{eq4}), and  we get a  local solution on the right of $r_o$.   The situation $w(r_o) >0$  cannot occur, since  this would imply that $\lim_{r\rightarrow r_o, r> r_o} {d\over dr} (|w^\prime |^\alpha w^\prime(r)) <0$ and $w^\prime $ coud not be  $\leq 0$ on the left of $r_o$ and $=0$ on $r_o$.

Proof of proposition \ref{propk}

We prove the result   for  equation (\ref{eq2}), with $k_o=1$ and $r_o=0$,  the changes to bring in the other cases are  given shortly at the end of the proof. 

The equation can also be written as 
$$ \left\{ \begin{array}{lc}
{d\over dr} (r^{(N-1)(1+\alpha)} |k^\prime |^\alpha k^\prime )(r)= -{(\alpha+1)r^{(N-1)(1+\alpha)} |k|^\alpha k(r)\over a}& {\rm in} \ \R^+\\
  k(0) = 1,  k^\prime (0) = 0.  & (4.7)
  \end{array}\right.$$
  or equivalently,   defining $\varphi_{p^\prime } (u) = |u|^{p^\prime -2} u $ and $p^\prime = {\alpha+2\over \alpha+1}$ as :  
   \begin{equation}\label{eqlap}
   k(r) = 1-\int_0^r \varphi_{p^\prime} \left({\alpha+1\over a s^{(N-1)(1+\alpha)} }\int_0^s t^{(N-1)(1+\alpha)} |k|^\alpha k (t) \ dt\right)  ds.
   \end{equation}
    We use the properties of the operator 
    \begin{equation}\label{Tlap}
     T(k) (r)= 1-\int_0^r \varphi_{p^\prime} \left({\alpha+1\over a s^{(N-1)(1+\alpha)} }\int_0^s t^{(N-1)(1+\alpha)} |k|^\alpha k (t) \ dt \right) d s
     \end{equation}
     which  satisfies  on $[0, \delta]$
     
     $$||T(k)-1||_\infty \leq\delta \left\vert \varphi_{p^\prime} \left({(\alpha+1)\delta ||u||_\infty ^{\alpha+1}\over a( (N-1) (1+\alpha)+1)} \right)\right\vert \leq c_1 \delta ^{p^\prime} ||u||_\infty \leq c_1 \delta ^{p^\prime} (||u-1||_\infty+1) $$
where $c_1 =   \left({(\alpha+1)\over a ((N-1) (1+\alpha)+1)} \right) ^{p^\prime-1}$

   If $\delta <\left ({1\over  3^{|\alpha|+1}c_1}\right)^{1\over p^\prime}$,   $T$ sends the ball  $\{ u\in {\cal C} ([0, \delta]),  ||u-1||_{{\cal C} ([0, \delta])}\leq {1\over 2}\}$ into itself.  We now prove that it is contracting. 
     We  observe that  for $k$ with values in $[{1\over 2}, {3\over 2}]$
   \begin{eqnarray*}
    {(\alpha+1)\over a ((N-1) (1+\alpha)+1)} \left( {1\over 2}\right)^{\alpha+1} \ s&\leq& 
 {\alpha+1\over a s^{(N-1)(1+\alpha)} }\int_0^s t^{(N-1)(1+\alpha)} |k|^\alpha k (t) \ dt\\
 & \leq& {(\alpha+1)\over a ((N-1) (1+\alpha)+1)}  \left({3\over 2}\right)^{\alpha+1}\ s, 
 \end{eqnarray*}
 
 and then by the mean value theorem  for $(u,v)\in B_{{\cal C} ([0, \delta])} (1, {1\over 2})$
 \begin{eqnarray*}
 \left\vert  \varphi_p^\prime  \left({\alpha+1\over a s^{(N-1)(1+\alpha)} }\int_0^s t^{(N-1)(1+\alpha)} u^{1+\alpha}(t)dt \right)\right. &-&\left.\varphi_{p^\prime} \left({\alpha+1\over a s^{(N-1)(1+\alpha)} }\int_0^s t^{(N-1)(1+\alpha)} v^{1+\alpha}  (t) dt\right)\right\vert \\
 &\leq& c_1 s^{p^\prime-1} |u^{\alpha +1}-v^{\alpha +1}  |_{L^\infty ([0, s])} \sup \left(\left({3\over 2}\right)^{-\alpha} , \left({1\over 2}\right)^{-\alpha}\right)\\
 &\leq &c_1 s^{p^\prime-1} |u-v| _{L^\infty ([0, s])}   \sup\left(\left( {3\over 2}\right)^{-\alpha} , \left({1\over 2}\right)^{-\alpha}\right)\\
&& \sup \left(\left({3\over 2}\right)^{\alpha} , \left({1\over 2}\right)^{\alpha}\right)\\
  &\leq & c_1 s^{p^\prime-1} |u-v| _{L^\infty ([0, s])}  3^{|\alpha|}
  \end{eqnarray*}
   This implies that   
   $$|T(u)-T(v)|\leq  c_1{ \delta^{p^\prime}  \over p^\prime }  |u-v| _{L^\infty ([0, \delta ])} 3^{|\alpha|} 
   \leq {1\over 3}  |u-v| _{L^\infty ([0, \delta ])}$$
   Then the fixed point theorem implies that there exists a unique fixed point in ${\cal C} ([0, \delta])$.

     In the case of equation ( \ref{eq3}) one is lead to consider 
   $$  T(k)(r) = k_o-\int_{r_o}^{r}\varphi_{p^\prime} \left({\alpha+1\over a s^{N^+} }\int_{r_o}^s t^{N^+} |k|^\alpha k (t)\  dt \right) ds $$
   with $N^+ = {(N-1)(1+\alpha)A\over a}$.
   
   For equation (\ref{eq5})  we shall consider 
   
    $$  T(k)(r) = k_o-\int_{r_o}^r\varphi_{p^\prime} \left({\alpha+1\over A s^{N^-} }\int_{r_o}^s t^{N^-} |k|^\alpha k (t)\  dt \right)ds$$
   with $N^- = {(N-1)(1+\alpha)a\over A}$.
    Finally for equation  (\ref{eq4})
        $$ T(k) (r)= k_o-\int_{r_o}^r\varphi_{p^\prime} \left({\alpha+1\over A s^{(N-1)(1+\alpha)} }\int_{r_o}^s t^{(N-1)(1+\alpha)} |k|^\alpha k (t)\  dt \right) ds.$$

   \bigskip

Proof of proposition \ref{propmM}

       We prove the local existence  by  proving  that for each $(w_o,w_o^\prime)$ with $w_o^\prime \neq 0$  and for all $ r_o>0$, there exists a neighborhood around $ r_o$  and a solution to the equation  which satisfies  the condition $(w( r_o), w^\prime ( r_o) )= (w_o,w_o^\prime)$.
      We suppose that $w_o^\prime \neq 0$ and we introduce the function 
      $$f_2(r, y_1, y_2) = M\left(-{m(y_2) (N-1)\over r }- {|y_1|^\alpha y_1 \over |y_2|^\alpha}\right)$$
       where $M$ and $m$ are respectively the functions 
       $$ M (x)= \left\{\begin{array}{c}
        {x\over A}\ {\rm if} \ x>0\\
         {x\over a} \ {\rm if } \ x<0
         \end{array}\right.$$
          and 
           $$ m (x)= \left\{\begin{array}{c}
        { A x }\ {\rm if} \ x>0\\
         {a x } \ {\rm if }\  x<0
         \end{array}\right.$$
    The functions $M$ and $m$ are lipschitzian, hence $f_2$ is lipschitzian with respect to $y = (y_1, y_2)$ around  $(w_o, w_o^\prime)$ when $w_o^\prime \neq 0$. 
     Let $f_1(r, y_1, y_2) = y_2$, and $f(y_1, y_2) = (f_1(y_1, y_2), f_2(y_1, y_2))$, 
     then the standard theory of ordinary differential equations implies that 
     $$\left\{ \begin{array}{c}
     (y_1^\prime, y_2^\prime) = f(y_1, y_2)\\
      (y_1, y_2)(r_o) = (w_o,w_o^\prime)
      \end{array}\right.$$
     has a unique solution around $(w_o,w_o^\prime)$ when $w_o^\prime\neq 0$.  Then $w = y_1$ is a  local solution of 
     \begin{equation} \label{eqmM}w^{\prime\prime } = M\left(-{m(w^\prime ) (N-1)\over r}-{|w|^\alpha w \over |w^\prime |^\alpha}\right)
     \end{equation}
       with the initial condition $w ( r_o) = w_o ,$ $w^\prime ( r_o) = w^\prime_o $. 
       
        If $w$ is a solution on $]r_1, r_2[$ and $\lim_{r\rightarrow r_2, r<r_2}w^\prime(r)$ exists and is $\neq 0$,  $w^{\prime\prime } $ has also a finite limit  from  the equation,  then $\lim_{r\rightarrow r_2, r< r_2} (y_1^\prime, y_2^\prime )(r) $ exists and is finite and one can  continue,  replacing $r_o$ by $r_2$ and $(w_o, w^\prime_o)$ by $(w(r_2), \lim_{r\rightarrow r_2, r<r_2}w^\prime(r))$.    If  $\lim_{r\rightarrow r_2, r<r_2}w^\prime(r))$  is zero, one gets $\lim_{r\rightarrow r_2} w^{\prime\prime } (r)= \pm \infty$ and then one cannot get a continuation, since the solutions of $(y_1^\prime, y_2^\prime) = f(y_1, y_2)$ must be ${\cal C}^1$.

   We prove the last facts concerning ${d\over dr} (|w^\prime |^\alpha w^\prime)$. Suppose that $w(r_2) >0$ and assume by contradiction that $\lim_{r\rightarrow r_2} {d\over dr} (|w^\prime |^\alpha w^\prime(r))\geq 0$. Then the equation on the left of $r_2$ is, since it is clear from the equation that $w^\prime$ cannot be nonnegative  : 
   $$A {d\over dr}\left({|w^\prime|^\alpha w^\prime \over 1+ \alpha}\right) + {a(N-1) \over r} |w^\prime |^\alpha w^\prime = -|w|^\alpha w $$
    which yields  a contradiction when $r\rightarrow r_2$. 
       
     Suppose now that $w(r_2)<0$ and  $\lim_{r\rightarrow r_2} {d\over dr} (|w^\prime |^\alpha w^\prime(r))\leq 0$, then   from the equation $w^\prime $ cannot be $\geq 0$ on the left of  $r_2$, and one is lead to solve   on the left  of $r_2$  : 
        $$a {d\over dr} \left({|w^\prime|^\alpha w^\prime \over 1+ \alpha} \right)+ {A(N-1) \over r} |w^\prime |^\alpha w^\prime = -|w|^\alpha w .$$
        This is absurd by passing to the limit when $r\rightarrow r_2$.

      Suppose that $w(r_1)>0$ and assume  by contradiction that      $\lim_{r\rightarrow r_1, r> r_1} {d\over dr} ( |w^\prime|^\alpha w^\prime )(r) \geq 0$, by the equation $w^\prime$  cannot be $\geq 0$ then this  equation is  on the right of $r_1$
      
         $$A {d\over dr}\left( {|w^\prime|^\alpha w^\prime \over 1+ \alpha} \right)+ {a(N-1) \over r} |w^\prime |^\alpha w^\prime = -|w|^\alpha w. $$         
      This is absurd by passing to the limit  when $r\rightarrow r_1$.   
      
         Suppose that $w(r_1)<0$ and that  
        $\lim_{r\rightarrow r_1, r>  r_1} {d\over dr} ( |w^\prime|^\alpha w^\prime )(r) \leq 0$, then the equation on the right of $r_1$ is 
            $$a {d\over dr} \left({|w^\prime|^\alpha w^\prime \over 1+ \alpha} \right)+ {A(N-1) \over r} |w^\prime |^\alpha w^\prime = -|w|^\alpha w . $$ 
     This is absurd  letting $r$ go to $r_1$.

      \bigskip
      
     Proof of proposition \ref{propr1}
     
We can observe that $ |k^\prime|^\alpha  k^\prime$ is differentiable for $r>0$  and has a limit $<0$ for $r\rightarrow 0$. Moreover we shall give some constant $\delta_1$  which depends only on $a$, $A$, $\alpha$, $N$  such that $k^\prime \neq 0$ and   ${d\over dr} \left(|k^\prime |^\alpha k^\prime\right)$ remains $<0$ on $]0, \delta_1[$. 

We begin to prove that   ${d\over dr} ( |k^\prime|^\alpha k^\prime)< 0$ around zero. One has  for $r>0$ 
$$(|k^\prime |^\alpha k^\prime)(r) = -{1+\alpha\over a r^{N_o}}\int_0^r (|k|^\alpha k)(s) s^{N_o} ds$$
where $N_o = (N-1)(1+\alpha)$, and then  $ ( |k^\prime|^\alpha k^\prime)$ is continuously differentiable for $r\neq 0$, as the primitive of some continuous function,  and 
$${d\over dr} (|k^\prime |^\alpha k^\prime ) (r)= {N_o(1+\alpha)\over a r^{N_o+1}} \int_0^r (|k|^\alpha k)(s) s^{N_o} ds -
{1+\alpha \over a}( |k|^\alpha k)(r). $$

For the point $0$, one has 
$$\lim_{r\rightarrow 0} {|k^\prime|^\alpha k^\prime(r)\over r} = -\lim_{r\rightarrow 0} {1+\alpha\over ar^{N_o+1}} \int_0^r  |k|^\alpha k (s)s^{N_o} ds  = -{1+\alpha\over a (N_o+1)}<0$$

 Using the fact that $k$ tends to $1$  when $r$ goes to zero we get that

$$\lim_{r\rightarrow 0}{d\over dr} (|k^\prime |^\alpha k^\prime )(r)  = {1+\alpha \over a} ({N_o\over N_o+1} -1)= -{(1+\alpha) \over A(N_o+1)}<0. $$
and then $|k^\prime |^\alpha k^\prime $ is ${\cal C}^1$ on $0$. 

 Moreover we prove that there exists a neighborhood  on the right of zero which depends only on the data, such that ${d\over dr} ( |k^\prime|^\alpha k^\prime)<0$ on it.  For that aim we begin to establish  some Lipschitz estimate on the solution with  some constant which depends only on  the data.   
  
  We have chosen $\delta$ (which depends only on $a$, $A$, $\alpha$, and $N$)  such that  for  $r\in [0, \delta]$, $k(r)\in [{1\over 2}, 
{3\over 2} ]$. We now observe  that $k^\prime$ is then bounded by 
  $$|k^\prime |^{\alpha+1}(r) \leq {1+\alpha \over a(N_o+1)} \left({3\over 2}\right)^{\alpha+1} r$$

   We have obtained that  there exists some constant $c_2$ which depends only on the constant $a$, N$, A$ such that 
   $|k^\prime |\leq c_2$
   on $]0, \delta[$. 
    We derive from  this that  on $[0, \delta]$
    $$|k(r)-1|\leq c_2  r, $$
    and also  that 
    $$|(|k|^\alpha k)(r) -1|\leq (1+\alpha) \sup (\left({3\over 2}\right)^\alpha, \left({1\over 2}\right)^\alpha ) c_2 r = c_3 r, $$
     and then 
    \begin{eqnarray*}
     |{d\over dr} (|k^\prime |^\alpha k^\prime )(r) + {1+\alpha \over a(N_o+1)} |&\leq& 
    {(1+\alpha )N_o\over a r^{N_o+1}} \int_0^r ||k|^\alpha k-1| (s)s^{N_o} ds \\
    &+&
{1+\alpha \over a}|( |k|^\alpha k)(r)-1| \\
&\leq & { c_3(1+\alpha)r\over a} \left({N_o\over N_o+2} +1\right)
\end{eqnarray*}

 We have obtained  that as long as $r< {N_o+2\over 2(N_o+1)^2 c_3}\equiv r_1$,  ${d\over dr} ( |k^\prime|^\alpha k^\prime)$ 
 remains negative (and then so does $k^\prime$). 
 This ends the proof of proposition \ref{propr1}. 
 \bigskip

 To finish the proof of proposition  \ref{propex},  i.e. to prove global existence's result, 
  suppose that $w$ is a solution on $[0, r_1[$. If $w^\prime (r_1)\neq 0$ we use proposition \ref{propmM}, if $w^\prime (r_1)=0$, using $\lim_{r\rightarrow r_1, \ r> r_1} {d\over dr} (|w^\prime |^\alpha w^\prime ) (r)  w(r_1) <0$ we consider  on the right of $r_1$, equation ( \ref{eq2})   if $w(r_1)>0$,  and equation ( \ref{eq4}) if $w(r_1)<0$. We have obtained a solution on $\R^+$. 
     
\bigskip
 We now prove that the solution $w$ is oscillatory :

\begin{prop}\label{proposc}

The solution of  ( \ref{eq1}) is oscillatory, ie, for all $r>0$ there exists $\tau> r$ such that $w(\tau) = 0$.
 
  \end{prop}
  
   Proof of proposition  \ref{proposc} :

   {\bf First step}

   We suppose that  $a= A$.  
      We follow the arguments in \cite{DM}. 
   
   We assume by contradiction that there exists $r_o$ such that $ w$ does not vanish on $[r_o, \infty[$. Then one can consider the function

$$y(r) = r^{(N-1)(1+\alpha)} {|w^\prime |^\alpha w^\prime(r) \over |w|^\alpha w(r)}, $$

which  satisfies the equation
 $$ y^\prime(r) = -{(\alpha+1)  r^{(N-1)(1+\alpha)}\over a}-{(\alpha+1)|y|^{\alpha+2}(r)  \over r^{(N-1)(1+\alpha)^2}} . $$
  Integrating between $r_0$ and $t$ one gets  that 
  $$ y(t) + (\alpha+1)\int_{r_0}^t {|y|^{\alpha+2}(r)  \over  r^{(N-1)(1+\alpha)^2} }dr = -{(\alpha+1)t^{(N-1)(1+\alpha)+1} \over a \left((N-1)(1+\alpha)+1\right)} + y(r_0)+{(\alpha+1)r_o^{(N-1)(1+\alpha)+1} \over a\left((N-1)(1+\alpha)+1\right)} . $$
   In particular we obtain that $y(t)\leq 0$  for $t$ large enough. 
   
     For the next step  it will be useful to remark that  if, in place of the equation,  we had  
   the inequation 
   $${d\over dr}\left( r^{(N-1)(1+\alpha)} |w^\prime |^\alpha w^\prime(r) \right)\leq { -r^{(N-1)(1+\alpha)} |w|^\alpha w\over a}, $$
   the conclusion would be the same. 
  \bigskip
  
    We obtain that 
    $ -y(t) = |y(t)|\geq C t^{(N-1)(1+\alpha) +1}$ for some constant $C>0$, 
    as soon as $t$ is large enough. 
    Let $k(t) = \int_{r_0}^t {|y|^{\alpha+2} (r)\over  r^{(N-1)(1+\alpha)^2} }dr $,  
     then using the previous considerations 
     $k(t) \geq c_1t^{N(1+\alpha)+2}$ for some positive constant $c_1$.
     
     Coming back to the equation, always for $t$ large  
    $$  (\alpha+1)k(t) \leq  |y(t)| = (k^\prime(t) t^{(N-1)(1+\alpha)^2})^{1\over \alpha+2} $$ and then, 
     $$(1+\alpha)^{\alpha+2} k^{\alpha+2}(t) \leq k^\prime(t) t^{(N-1)(1+\alpha)^2}.$$
     Integrating between $t$ and $s$, $s> t$, we obtain  that for some positive constant $c_2$ 
     $${1\over k^{\alpha+1} (t)} -{1\over k^{\alpha+1} (s) } \geq c_2\left( {1\over t^{(N-1)(1+\alpha)^2-1}} -{1\over s^{(N-1)(1+\alpha)^2-1}} \right).$$
     Letting $s$ go to infinity 
      $${1\over k^{\alpha+1} }(t) \geq c_2{1\over t^{(N-1)(1+\alpha)^2-1}} .$$
      From this  one gets a contradiction with $k(t) \geq c_1t^{N(1+\alpha)+2}$. This ends the proof of the first step.

     \bigskip 
     {\bf Second step :}
      
      $a < A$.

      We argue on the model  of  \cite{BEQ}.

      We suppose as in the first step that there exists $r_o$ such that $w$ does not vanish on $[r_o, \infty[$.  
    
       We begin to prove that if $w >0$ for $r\geq r_o$, then  for $r\geq r_o$ 
       $${d\over dr} \left(r^{(N-1)(1+\alpha)} |w^\prime |^\alpha w^\prime(r)\right) \leq { -r^{(N-1)(1+\alpha)} |w|^\alpha w(r)\over a}, $$
        and then following the previous arguments  in the first step we obtain that if $y(r) = r^{(N-1)(1+\alpha)} {|w^\prime |^\alpha w^\prime(r) \over |w|^\alpha w(r)} $ then 
         $$ y(t) + (\alpha+1)\int_{r_o}^t {|y|^{\alpha+2} (r)\over  r^{(N-1)(1+\alpha)^2} }dr \leq  -{t^{(N-1)(1+\alpha)+1} \over  a\left((N-1)(1+\alpha)+1\right)} + y(r_o)+{r_o^{(N-1)(1+\alpha)+1} \over  a\left((N-1)(1+\alpha)+1\right)} , $$
         a contradiction if $y>0$ for $t$ large enough. 
         
          To prove  that  $${d\over dr}\left( r^{(N-1)(1+\alpha)} |w^\prime |^\alpha w^\prime(r)\right) \leq { -r^{N-1)(1+\alpha)} |w|^\alpha w(r)\over a}, $$
         let us note that in the case $w^\prime\leq 0$ and ${d\over dr} ({|w^\prime|^\alpha w^\prime \over 1+\alpha})\leq 0$ equality holds in the previous inequality, if $w^\prime\geq 0$ and ${d\over dr} ({|w^\prime|^\alpha w^\prime \over 1+\alpha}) \geq 0$ the equation is impossible . 
          For the other cases,  
         we assume first that $w^\prime \leq 0$,  this implies  if ${d\over dr} \left({|w^\prime|^\alpha w^\prime \over 1+\alpha}\right)\geq 0$ that 
           \begin{eqnarray*}
           a{d\over dr} \left({|w^\prime|^\alpha w^\prime \over 1+\alpha} \right)+ { a (N-1)\over r} |w^\prime |^\alpha w^\prime &\leq & \left(A {d\over dr} ({|w^\prime|^\alpha w^\prime \over 1+\alpha}) + { a (N-1)\over r} |w^\prime |^\alpha w^\prime \right) \\
          &\leq & -|w|^\alpha w, 
           \end{eqnarray*}
           which  implies the result. 
           
            If $w^\prime \geq 0$ and ${d\over dr} ({|w^\prime|^\alpha w^\prime \over 1+\alpha})\leq 0$
                \begin{eqnarray*}
           a{d\over dr} \left({|w^\prime|^\alpha w^\prime \over 1+\alpha} \right)+ { a (N-1)\over r} |w^\prime |^\alpha w^\prime &\leq & a {d\over dr} ({|w^\prime|^\alpha w^\prime \over 1+\alpha} )+ { A (N-1)\over r} |w^\prime |^\alpha w^\prime  \\
                    & \leq & -|w|^\alpha w.
           \end{eqnarray*}
           This also implies the result.

             We  now  assume that $w<0$ on $[r_o, \infty[$. Then we prove that there exists $r^\star$ such that $w^\prime (r^\star) = 0$ and $w^\prime >0$ on $]r^\star, \infty[$. 
             
              Indeed by the equation if $w^\prime (r^\star) = 0$, by proposition \ref{propmM} $\lim_{r\rightarrow r^\star, r>  r^\star} {d\over dr} (|w^\prime|^\alpha w^\prime ) (r) >0$. This implies that $w^\prime$ is  increasing on $r^\star$, then $w^\prime$ is $>0$ on a neighborhood on the right of $r^\star$. 
              
            If there exists $r^\prime > r^\star$ such that 
               $w^\prime (r^\prime ) = 0$, we argue as before and then $w^\prime >0$ after $r^\prime$. 
               
                From these remarks, it is sufficient to  discard  $w^\prime <0$ on $[r_o, \infty[$. Then in that case necessarily ${d\over dr} (|w^\prime |^\alpha w^\prime)  >0$ on $[r_o, \infty[$ by the equation, and then $w$ satisfies 
                  ${d\over dr} (|w^\prime |^\alpha w^\prime(r) r^{N^-}) = -(1+\alpha) {r^{N^-} |w|^\alpha w(r)\over a} >0$. 
                    Let $g(r) \equiv  (|w^\prime |^\alpha w^\prime(r) r^{N^-}) $,  $g$ is monotone increasing , and since $w^\prime <0$,  it has a limit $c_1\leq 0$ at $+\infty$. On the other hand, since $w^\prime <0$ there exists $c_2\in [-\infty, 0[$ such that $\lim_{r\rightarrow +\infty} w(r) = c_2$, then from the equation satisfied by  $w$, 
                    $\lim_{r\rightarrow +\infty} g^\prime (r) = +\infty$,  which is a contradiction with $\lim_{r\rightarrow +\infty} g(r) = c_1 \leq 0$.

                    Finally $w^\prime >0$ after $r_o$. 
                  
                    We recall that    $N^+ = {A(1+\alpha)  (N-1)\over a} $. 
                       Distinguishing the cases ${d\over dr} ( |w^\prime|^\alpha w^\prime)>0$ and ${d\over dr} ( |w^\prime|^\alpha w^\prime) <0$  on the right of $r_o$  and arguing as we already did before we obtain that 
                     $w$     satisfies        $${d\over dr} (|w^\prime |^\alpha w^\prime(r) r^{N^+}) \leq -{(1+\alpha) r^{N^+} |w|^\alpha w(r) \over a}$$
                     
                   Then defining  
                    $$y(r) = r^{N^+}{ |w^\prime |^\alpha w^\prime(r) \over |w|^\alpha w(r)} $$
                     one has 
                    \begin{equation}\label{eqb}y^\prime(t) + {(\alpha+1)|y(r)|^{\alpha+2} \over r^{(N^+)(\alpha+1)}} + {(\alpha+1)r^{N^+}\over a} \leq 0.
                    \end{equation}

                     Hence integrating between $r_o$ and $t$ one gets  for some constant $c_1>0$ 
                     $$|y(t)|= -y(t) \geq c_1  t^{N^++1}. $$

                     Let $k(t) = \int_{r_o}^t {|y|^{\alpha+2}(r) \over r^{N^+(\alpha+1)}} dr \geq c t^{N^++ \alpha+3}$. 
                     From the  equation (\ref{eqb}) integrated between $r_o$ and $t$ ,  using 
                     $$k^\prime (t) = {|y|^{\alpha+2}(t) \over t^{N^+(\alpha+1)}}, $$
                      we get 
                      $$(\alpha+1)^{\alpha+2}k^{\alpha+2}(t) \leq k^\prime(t) t^{N^+(\alpha+1)}, $$
                      hence   for some positive constant $c_2$
                      $$k^{-(\alpha+1)}(t) -k^{-(\alpha+1)} (s)\geq c_2( t^{-N^+(\alpha+1)+1}-s^{-N^+(\alpha+1)+1})$$
                      for $s>t$. 
                      Letting $s$ go to infinity and using $\lim k(t) = +\infty $,  one derives that 
                       $$k^{-(\alpha+1)}(t) \geq c_2t^{-N^+(\alpha+1)+1}, $$
                        which is a contradiction with $k(t) \geq c_1t^{N^++ \alpha+3}$. 
                       
                        We have obtained that $w$ is oscillatory.  This ends the proof of  proposition \ref{proposc}.
                        
                         For the sake of completeness, we give some property of the function $w$ inherited from the property of the eigenfunctions in the viscosity sense \cite{BDr} :

\begin{lemme}
Between two successive   zeros of $w$, there exists a unique zero of $w^\prime$. 
\end{lemme}

 Proof

  Suppose that $w$ is constant sign on $B(0,t)\setminus\overline{B(0, s)}$, $s<t$ and $w(s)=  w(t)=0$, then $w_1 (x) = w(\mu^{1\over 2+\alpha} x)$ is an eigenfunction for  one of the first demi-eigenvalue $\mu = \lambda^+(B(0,t)\setminus\overline{B(0, s)})$  if $w>0$,  or $ \mu = \lambda^-(B(0,t)\setminus\overline{B(0, s)})$   if $w<0$. Then by the uniqueness of the first eigenfunction  in the radial case,   if $w>0$, by remark \ref{remaru}, $w$ is increasing on $[s, r_w]$ and decreasing on $[r_w,     t]$ and $r_w$ is the unique point on which $w^\prime =0$.   
  We argue in the same manner when  $w<0$,  using the fact that in that case  $w$ is decreasing on $[s, r_w]$ and increasing on $[r_w, t]$.

                        \bigskip
                         In the sequel we shall  denote by $w^+$  the  radial solution given by proposition \ref{propex}  of 
                         $$\left\{ \begin{array}{cc}
                          |w^\prime |^\alpha {\cal M}_{a, A} (r, w^\prime, w^") = -|w|^\alpha w& \\
                           w(0) = 1,\ w^\prime (0) = 0.& 
                           \end{array}\right.$$ 
                          
                       And we denote by $w^-$ the radial solution of 
                         $$\left\{ \begin{array}{cc}
                          |w^\prime |^\alpha {\cal M}_{a, A} (r, w^\prime, w^")  = -|w|^\alpha w&\  \\
                           w(0) = -1,\ w^\prime (0) = 0.& 
                           \end{array}\right.$$ 
                            The  proof of the existence  and uniqueness of $w^-$ is obtained by the same arguments used for   $w^+$. 
                           The results  in proposition \ref{proposc} can be adapted to the case of $w^-$, and then   we also get that $w^-$ is oscillatory. 
                          \section{Eigenvalues and eigenfunctions}
                           In this  section we  prove the existence of an infinite numerable set of eigenvalues  for the radial operator defined in  equation (\ref{eqpucci}).  These eigenvalues are simple and isolated. We begin with some  properties of  the eigenfuntions. 
                           
                            \begin{prop}\label{propuprime0}
  
   Suppose that $u$ is a  radial  viscosity solution of 
    $$\left\{ \begin{array}{lc}
    \tilde F(r,  u^\prime , u^") = -\mu |u|^\alpha u& {\rm in} \ B(0,1)\\
    u(1) = 0, u(0) >0. & \ 
    \end{array}\right.$$
    
 Then $0$ is a local maximum for $u$,   $u$ is ${\cal C}^2$  on a neighborhood $]0, r_o[$ of zero, is ${\cal C}^1$ on $[0, r_o]$ and  $u^\prime (0) =  0$.
     \end{prop}
     
      Proof  of proposition \ref{propuprime0}
      
      First let us note that $\mu >0$, because if not the maximum principle would imply that $u\leq 0$. 
      
       Since $u$ is continuous  there exists some neighborhood $B(0, r_o)$ on which 
      $$\tilde F(r,  u^\prime , u^") <0$$
      Then   using  the   comparison principle   for such operators,   and remarking that  positive constants  are sub-solutions,  one gets that $u(r) \geq u(r_1)$ on $B(0, r_1)$,  if $r_1< r_o$. This implies in particular that $u$ is decreasing from zero, and $0$ is a local maximum.  We   now prove that $u$ is ${\cal C}^1$ around zero and ${\cal C}^2$ on a neighborhood of $0$, except on $0$.  
   
          Let $r_1$ be the first zero of $u$. Then $u>0$ on $B(0, r_1)$ and $\lambda^+(B(0,  r_1) )= \mu, $ by proposition \ref{mupropre}. 
      Let $w^+$ be the ${\cal C}^1$ solution in  proposition \ref{propex} and $\beta^+_1$ its first zero,  (it exists according to proposition \ref{proposc}). Define  
        $$ v(r) = w^+({\beta_1^+ r\over  r_1}).$$
        Then $v>0$ on $B(0, r_1)$ and $v$ is an eigenfunction  in $B(0,  r_1)$ for the eigenvalue $\left({\beta_1^+\over  r_1}\right)^{2+\alpha}$, in particular $\lambda^+(B(0,  r_1)) =\mu =  \left({\beta_1^+\over  r_1}\right)^{2+\alpha}$, and by the uniquenes of the first radial  eigenfunction $ >0$  in proposition \ref{propun}, there exists some constant $c>0$ such that 
        $ u = c v$  on $B(0, r_1)$. 
        In particular $u$ is ${\cal C}^2$ on each point where $u^\prime$ is different from zero and  ${\cal C}^1$ everywhere on $B(0,  r_1)$.     This proves in particular,  since $u$ is ${\cal C}^1$ on $B(0,  r_1)$ and $u$ has a maximum on $0$, that $u^\prime (0)=0$. 
         
Of course the symmetric result  holds  for $u$ such that $u(0) <0$. 

\bigskip
We now present an improvement of proposition \ref{propinc} which will be used in the proof of corollary \ref{cormuk}

   \begin{prop}\label{propts}
   Suppose that $s<t<1$
   
    Suppose that there exists some  eigenfunctions for the annulus $B(0, 1) \setminus\overline{B(0,s)}$ and for $B(0,1)\setminus \overline{B(0, t)}$,  which are ${\cal C}^2$ on each point where their  first derivative is different from $0$,  and ${\cal C}^1$ anywhere,  then 
    $\lambda^\pm (B(0, 1) \setminus\overline{B(0,s)})< \lambda^\pm (B(0, 1) \setminus\overline{B(0,t)}) $. 
    \end{prop}
    
    Proof 
    
    Suppose by contradiction that  $\lambda^\pm (B(0, 1) \setminus\overline{B(0,s)})= \lambda^\pm (B(0, 1) \setminus\overline{B(0,t)}) $. that we shall denote for simplicity by $\lambda^\pm$. 
    Let $\varphi$ and $u$ be solutions  of 
    the equation 
    $$ \tilde F(r,  \varphi^\prime , \varphi^")+  \lambda^\pm |\varphi|^\alpha \varphi = 0$$
     which are ${\cal C}^2$ on each point where their  first derivative is different from $0$,  and ${\cal C}^1$ anywhere, 
     with $\varphi = 0$ on $\{r=1\}$ and $\{r= s\}$, and $u=0$ on $\{r=1\}$ and $\{r=t\}$. 
     To fix the ideas we also assume that $\varphi$ and $u$ are positive (and then we   replace  $\lambda^\pm $ by   $\lambda^+$) . 
     
      Using the same arguments as in propositions \ref{propk} and \ref{propmM},  since  $\varphi (1) = u(1)=0$ and $u^\prime (1) <0$, $\varphi^\prime (1)<0$,  by  uniqueness  there exists some constant $c>0$  such that  $\varphi = cu$ as long as $\varphi^\prime$  or $u^\prime$ is  different from zero.  By remark \ref{remaru} there exists exactly one point  
      $r_u$ on $]t,1[$ for which   $u^\prime (r_u)=0$ and it is a    global strict maximum for $u$ on $]t,1[$. By uniqueness, $\varphi^\prime (r_u)=0$ and $r_u$ must also be  a global strict maximum for  $\varphi$ on $]t,1[$. 
      Then the equation  satisfied by $u$ and $\varphi$ on the left of $r_u$, is equation (\ref{eq3}). By  local uniqueness   of solutions to (\ref{eq3}) one gets that $u = c \varphi$ on the left of $r_u$ and this is true as long as $u^\prime$ or $\varphi^\prime $ is different from $0$, hence  at least on $]t,1[$. We get a contradiction since $u = 0$ on $\{r=t\}$ and $\varphi(t)\neq 0$.

           \bigskip
           
            We  now prove the existence of  a numerable set  of eigenvalues.

                          The result in proposition \ref{proposc}  implies that there exists a sequence $\beta_k^\pm$ of increasing sequence of  zeros of $w^\pm$.  
                          
                                                     We now consider 
                           $u_k^\pm(r) =  w^\pm(\beta_k^\pm r)$.
                           Then $u_k^\pm$  is an eigenfunction on $B(0,1)$ for the eigenvalue $\mu_k^\pm :=(\beta_k^\pm)^{\alpha+2}$
                           and it has $k-1$ zeros  inside the ball, say 
                           $r_i \equiv{\beta_i^\pm\over \beta_k^\pm}$, $i\in [1, k-1]$.   We need to prove that they are the only eigenvalues   :

                                       \begin{prop}\label{propri}
                                       
                                       The set of eigenvalues of the operator is  the set $\{\mu_k^\pm, \ k\geq 1\}$. 
                                       These eigenvalues are simple in the  following sense : 
                           Suppose that $v$ is some eigenfunction for the eigenvalue $\mu_k^\pm$,  which is   ${\cal C}^1$  and ${\cal C}^2$ on each point where the first derivative is different from $0$, then 
                         there exists some constant $c>0$ such that 
                            $v =c w^\pm((\mu_k^\pm)^{1\over 2+\alpha} \cdot)$.
                            
                                                         \end{prop}
                             
                             Proof of proposition  \ref{propri}
                             
                             Let $\mu$ be an eigenvalue.
                             Let $v$ be a corresponding eigenfunction, that we suppose  to fix the ideas such that $v(0)>0$. Necessarily since $v$ is radial and ${\cal C}^1$, $v^\prime (0) = 0$.  Let  $z(\cdot) = {v(\mu^{-1\over 2+\alpha}\cdot)\over v(0)}$. Then $z$ satisfies equation (\ref{eq2}) and  by uniqueness $z = w^+$ on $[0, \mu^{1\over 2+\alpha})$. This implies that $\mu^{1\over 2+\alpha}$ is one of the zeros of $w$.  This proves also the simplicity of the eigenvalue $\mu$. 
                             The fact that the eigenvalues are  isolated is a consequence  of the properties of the zeros of $w^+$.  
                             \bigskip

                             The following corollary   are not necessary  for the present paper, they will be useful    for the bifurcation results announced  in the final  concluding section  :

       \begin{cor}\label{cormuk}
       There is uniqueness (up to a  positive multiplicative constant) of the  $k$-th eigenfunction. As a consequence 
        one has $\mu_k^-< \mu_{k+1}^+$ and $\mu_k^+< \mu_{k+1}^-$. 
        \end{cor}

    Proof

    It is sufficient to prove that $\beta_k^+ < \beta_{k+1}^-$ and $\beta_k^-< \beta_{k+1}^+$. 
    
  We begin to prove that $\beta_1^+ < \beta_2^-$. One has 
    $$\lambda^-(]\beta_1^+, \beta_2^+[) = 1 = \lambda^- (]0, \beta_1^-[)=1 $$
     Suppose first that $ \beta_2^+< \beta_1^-$, then  it contradicts proposition \ref{propinc}
    
    If $\beta_2^+= \beta_1^-$,   one has  a contradiction with  proposition \ref{propts}. 
         
           We consider the case $k\geq 2$.  
      
       Suppose by contradiction that $\beta_k^- < \beta_{k+1}^+< \beta_{k+2}^+\leq \beta_{k+1}^-$, 
       and in a first time we assume that  $\beta_{k+2}^+< \beta_{k+1}^-$. In that case one would have 
        $$\lambda^\epsilon (]\beta_{k+1}^+, \beta_{k+2}^+[) = \lambda^\epsilon (]\beta_k^-, \beta_{k+1}^-[)=1, $$
         where $\epsilon =  sign  (-1)^{k+1}$, this would then contradict proposition \ref{propinc}.

        In a second time if we assume that $\beta_{k+2}^+ = \beta_{k+1}^-$, this contradicts proposition \ref{propts}

   In the same manner  we should  prove that $\beta_{k+1}^+ < \beta_{k+2}^-$.

  For the sake of completeness we  finish this section with  some additional property  of the eigenvalues . This result is an analogous of one  result   in \cite{BEQ}. 

\begin{prop}\label{propgap}
The gap between  the two first half eigenvalues is larger than between the second ones : 

 $${\mu_1^-\over \mu_1^+} \geq {\mu_2^-\over \mu_2^+}. $$
 \end{prop}
 {\rm Proof of proposition \ref{propgap}

 Let $\varphi_i^\pm$ $i = 1,2$ be the eigenfunctions associated with $\mu_i^\pm$ with $\varphi_i^\pm (0) = \pm 1$. 
 
  Let $r^+$ be the first zero of $\varphi_2^+$, $r^-$ the first zero of $\varphi_2^-$. We prove that 
  $ r^-\geq  r^+$. indeed, suppose by contradiction that $r^-<  r^+$, and define 
  $$A^+ = \{ r, \ r^+ < r< 1\}$$
  and 
   $$A^- = \{ r, \ r^- < r< 1\}$$
    then $A^+ \subset A^-$   and then 
    $$\lambda^-(A^+) = \mu_2^+ \geq   \lambda^-(A^-) > \lambda^+ (A^-) = \mu_2^-$$
     and 
      $$\lambda^+ (B_{r^+}) = \mu_2^+ < \lambda^+ (B_{r^-})< \lambda^- (B_{r^-}) = \mu_2^-.$$
       We have obtained a contradiction. 
       
       Moreover let us consider 
       $$\psi (x) = \varphi_2^+ (r^+  x)$$
        Then $\psi$ is a  radial solution on $B(0,1)$ of 
         $$ |\psi^\prime  |^\alpha {\cal M}_{a, A} (r,\psi^\prime , \psi^")= - (r^+)^{2+\alpha} \mu_2^+|\psi|^\alpha \psi, $$
          which implies since $\psi (1) = 0$,  that $(r^+)^{2+\alpha} \mu_2^+= \mu_1^+$,  by the definition of the first  half eigenvalue. 
          
           In the same manner 
          $$ (r^-)^{2+\alpha} \mu_2^-= \mu_1^-,$$
           and then 
           $${\mu_1^-\over \mu_2^-}  \geq {\mu_1^+\over \mu_2^+}, $$
           this yields the result. 
           
           \section{ The continuity of the spectrum  with respect to the parameters.}
           
            In this section we  let vary $\alpha\in ]-1, \infty[$ and $a\in [0, A]$  and for that reason we denote by $\tilde F_{\alpha,a}$ the operator $\tilde F$ defined before. We denote by $\mu_k^\pm (\alpha, a)$ the corresponding eigenvalues. 
     In order to prove the continuity of the map $(\alpha, a)\mapsto \mu_k^\pm(\alpha, a )$,  we  begin to  establish the boundedness of the eigenvalues $\mu_k^\pm (\alpha, a)$ when $\alpha $ belongs to some compact set of $]-1, \infty[$ and $a\in [0, A]$. 
     
   \begin{prop}\label{propbc}
   We suppose that $ a=A=1$. Let $\lambda_{eq, \alpha}  (]c,b[)$ be the first "radial"  eigenvalue for the set $B(0, b)\setminus \overline{B(0,c)}$ and for the operator $u\mapsto - {d\over dr} {|u^\prime |^\alpha u\over 1+\alpha} - {N-1\over r} u^\prime$. Then   there exists some continuous   function $\varphi(\alpha)$,  bounded on every compact set of  $[-1, \infty [$,  such that 
     $$\lambda_{eq, \alpha}(]c,b[)\leq \varphi(\alpha ) (b-c)^{-2-\alpha}.$$
     \end{prop}
     
      \begin{cor}\label{corbc}
       We assume that $a< A$. Then  
       $$\lambda^+_{a,A, \alpha} (]c,b[)\leq a\varphi(\alpha ) (b-c)^{-2-\alpha}.$$ 
           \end{cor}
           \begin{cor}\label{cormuk}
           
           For all $k\geq 1$
            $$\mu_k^+(\alpha, a) (B(0,1)) \leq a \varphi(\alpha) k^{2+\alpha},$$
            and 
            $$\mu_k^-(\alpha, a) (B(0,1)) \leq a \varphi(\alpha) (k+1)^{2+\alpha}.$$
            \end{cor}
           
           Proof of  proposition \ref{propbc}
         
         Let us note that one can  also use the following  result for general operators  satisfying the hypothesis in section  2,  proved in \cite{BD3}  : 
         {\it  There exists some constant $C$ which depends on $a$, $A$, $N$ such that  if $R$ is the radius of some ball included in  $\Omega$ then }
          $$\lambda^\pm (\Omega ) \leq {C\over R^{\alpha+2}}.$$
           But we shall  give a  more precise estimate here : 
           
            For the radial case,
         one  can easily see that $$\lambda_{eq, \alpha}= \inf_{u\in W_0^{1, 2+\alpha} (]c,b[)}{ \int_c^b |u^\prime |^{2+\alpha} (r) r^{(N-1)(1+\alpha)} dr \over  \int_c^b |u |^{2+\alpha} (r) r^{(N-1)(1+\alpha)} dr }.$$
             
Let us consider the function 
$u(r) = (r-c) (b-r)$.
 We need to  get an upper bound for 
 $$I = \int_c^b |2r-(c+b)|^{2+\alpha} r^{(N-1)(1+\alpha)} dr , $$
 and to get a lower bound for   
 $$ J = \int_c^b (r-c)^{2+\alpha} (b-r)^{2+\alpha} r^{(N-1)(1+\alpha)} dr . $$
 For the first  integral we use the inequality 
 $r^{(N-1)(1+\alpha)} \leq  2^{|1-(N-1)(1+\alpha)|}(r-\left({c+b\over 2}\right))^{(N-1)(1+\alpha)} + \left({c+b\over 2}\right)^{(N-1)(1+\alpha)}) $.
 
 In the following  $c(\alpha, N)$ is some constant which  can vary from one line to another but is bounded for $\alpha \in [-1, M]$. 
 We obtain that 
 \begin{eqnarray*}
 J &\leq& c(\alpha, N) \left(\int_c^b |r-\left({c+b\over 2}\right)|^{2+\alpha + (N-1)(1+\alpha)}dr\right.\\
 &+& \left.  \left({c+b\over 2}\right)^{(N-1)(1+\alpha)} \int_c^b|r-\left({c+b\over 2}\right)|^{2+\alpha} dr \right)\\
 &\leq &c(\alpha , N) \left((b-c)^{3+\alpha + (N-1)(1+\alpha)}+ (c+b)^{(N-1)(1+\alpha)}(b-c)^{3+\alpha}\right)
 \\
 &\leq & c(\alpha, N) (b-c)^{3+\alpha} (c+b)^{(N-1)(1+\alpha)}.
 \end{eqnarray*}
To minorize $I$  we use 
  $$ r^{(N-1)(1+\alpha)} \geq 2^{-|1-(N-1)(1+\alpha)|}\left((r-c)^{(N-1)(1+\alpha)} + c^{(N-1)(1+\alpha)}\right) $$
  and then 
  
  \begin{eqnarray*}
  I&\geq & c(\alpha, N) \int_c^b \left((r-c)^{2+\alpha + (N-1)(1+\alpha)} (b-r)^{2+\alpha} + c^{(N-1)(1+\alpha)} (r-c)^{2+\alpha}(b-r)^{2+\alpha}  dr\right) \\
  &\geq& c(\alpha, N) (b-c)^{5+ 2\alpha + (N-1)(1+\alpha)} B(
 N(1+\alpha)+2, 3+\alpha)\\
 &+ &c^{(N-1)(1+\alpha)} (b-c)^{5+ 2\alpha} B(3+\alpha, 3+\alpha)\\
  &\geq & c(\alpha, N) (b-c)^{5+2\alpha} b^{(N-1)(1+\alpha)}
  \end{eqnarray*}
   where in the previous lines, $B$ denotes the Euler function. We have obtained the result. 
    \bigskip

    Proof of corollary \ref{corbc}

     We use the inequality  in proposition \ref{proplam}
     $$\lambda^+ (B(0, b)\setminus \overline{B(0, c)}) \leq a \lambda_{eq} (B(0, b)\setminus \overline{B(0, c)}) $$

      Proof of  corollary \ref{cormuk}

      Let us recall that we have denoted by  $(r_i)_i$ the zeros of  the eigenfunction $ \varphi_k^+$. 
       $\mu_k^+ (B(0,1))$ coincides with  $\lambda^+ (B(0, r_1))$ and with $\lambda^+ (B(0, r_{i+1})\setminus \overline{B(0, r_i))}= \mu_1^+ (B(0, r_{i+1})\setminus \overline{B(0, r_i)})$,  for all $i \in [ 1, k]$. 
       Now, either $r_1\geq {1\over k}$,
    or    there exists $i_o\geq 2$ such that $r_{i_o+1}-r_{i_o} \geq {1\over k}$ . In each of the cases we get the result. 
         Concerning $\mu_k^-$ we use the inequality 
         $\mu_k^-\leq \mu_{k+1}^+$ in corollary \ref{cormuk}.

            \begin{prop}\label{contalpha}
            Let $M>0$ be given. 
             Suppose that $(\alpha_n, a_n) \rightarrow (\alpha, a)\in ]-1, M[\times [0, A] $, then 
             $\mu_k^\pm (\alpha_n,a_n)\rightarrow \mu_k^\pm (\alpha,a )$. 
             \end{prop}
             
              Proof  of proposition \ref{contalpha}
             
                     By corollary \ref{cormuk}, the sequence   $(\mu_k^\pm (\alpha_n, a_n))_n$ is bounded,  so we  can extract from it a subsequence,  denoted in the same manner for simplicity,  such that $\mu_k^\pm (\alpha_n,a_n) \rightarrow \mu$, for some $\mu\in \R^+$.    We fix the integer $k$. Let $\varphi_n$ be  such that $\varphi_n(0)=1$,   and  
                $$\left\{ \begin{array}{lc}
                \tilde F_{\alpha_n, a_n} (r, \varphi_n^\prime, \varphi_n^") + \mu_k^+ (\alpha_n, a_n) |\varphi_n|^{\alpha_n}\varphi_n= 0&\ {\rm in } \ B(0,1)\\
                 \varphi_n (1)= 0 & 
                 \end{array} \right.
           $$
           
           Using the compactness results in  corollary  \ref{comp} one can extract from $(\varphi_n)$ a subsequence  which will be denoted in the same manner for simplicity, which converges uniformly  to a  viscosity solution $\varphi$ of 
              $$\left\{ \begin{array}{lc}
                \tilde F_{\alpha,a } (r, \varphi^\prime, \varphi^") +\mu |\varphi|^{\alpha}\varphi= 0&\ {\rm in } \ B(0,1) \\
                 \varphi (1)= 0 & 
                 \end{array} \right.
           $$ 
           By the uniform convergence, $\varphi$ is not identically zero and $\varphi(0)=1$. Then $\mu$ is  some eigenvalue. We must prove first  that  that $\varphi$ has $k-1$ zeros,  secondly  that 
 $\varphi$ is ${\cal C}^1$ and ${\cal C}^2$ on every point where the first derivative is different from zero.           
                       Let   $j$ be such that  $(r_i)_{1\leq i \leq j-1}$ are  the zeros of $\varphi$. By  remark \ref{remhopf} in section 3,    $\varphi$ changes sign on each of them. As a consequence there exists $\delta >0$ such that  for all $i\in [1, j-1]$, on $[r_i-\delta, r_i+\delta]$,  $\varphi$ has no other zero than $r_i$ and on $[r_{i-1}+\delta, r_i-\delta]$ $\varphi $ has no zero. 
            From $\varphi (r_i-\delta) \varphi (r_i+ \delta) <0$,  one has for $n$ large enough 
           $\varphi_n (r_i-\delta) \varphi_n (r_i+ \delta) <0$, and then $\varphi_n$ has at least one zero in   $]r_i-\delta, r_i+\delta[$. 
            In the same manner there exists $m>0$ such that 
           $|\varphi|> m$ on  every $[r_{i-1}+\delta, r_i-\delta]$ , which implies by the uniform convergence of $\varphi_n$ towards $\varphi$ that $\varphi_n$ cannot have a zero in this intervall.  As a consequence $k\geq j$. 
           Moreover by the strict monotonicity  of $\varphi$ on $[r_i-\delta, r_i+\delta]$,  $\varphi_n$ is also monotone for $n$ large enough. This  implies in particular the uniqueness of zero of $\varphi_n$ on that intervall. Finally $j = k$. 
           
           There remains to prove that $\varphi$ is "regular", i.e. that $\varphi$ is ${\cal C}^2$ on each point where the first derivative is different from zero, and ${\cal C}^1$ anywhere.

           Suppose that $\bar r< \bar t$ are two successive  zeros of $\varphi$,  then for $n$ large enough,  there exists $r_n< t_n$ two successive  zeros of $\varphi_n$ which converge respectively to $\bar r,\bar t$. Moreover $\varphi_n$  (respectively $\varphi$) has constant   sign on  $]r_n, t_n[$ (respectively $]\bar r, \bar t[$).  One can assume without loss of generality that  this sign is negatif.

       We need to prove that $\varphi$ is  "regular " on $[\bar r,\bar t]$. 
             Let $r_n^\prime $ be the unique zero of $\varphi_n^\prime$ on $]r_n, t_n[$. 
              Then $\varphi_n$ is the unique fixed point  on $]r_n, r^\prime_n[$ ,  of the operator  
        $T_n$ defined as 
        $$T_n(w)(r)= \varphi_n(r_n^\prime )-\int_{r_n^\prime} ^r \varphi_{p^\prime} \left({(1+\alpha_n)\mu_k^+ (\alpha_n, a_n, A) \over A s^{N_n^-}} \int_{r_n^\prime }^s |w|^{\alpha_n} w(t) t^{N_n^-} dt\right) ds, $$
        where $N_n^- = {a_n(N-1)(1+\alpha_n)\over A}$.  
 One can prove as it is done in the proof of proposition \ref{propr1}  that there exists some neighborhood $]r^\prime _n-\delta , r_n^\prime [$  with $\delta$ which does not depend on $n$,  such that on the left of $r_n$,  $\varphi_n^\prime <0$ and $\varphi_n^{\prime\prime} >0$. 
 
  In the same manner  
         $\varphi_n$ is the unique fixed point of $T_n$  on   $]r_n^\prime, r_n^\prime+\delta[$  defined as 
          $$T_n(w) (r)= \varphi_n(r_n^\prime )-\int_{r_n^\prime}^r \varphi_{p^\prime} \left((1+\alpha_n){\mu_k^+ (\alpha_n, a_n, A) \over A s^{N_{0,n}} }\int_{r_n^\prime }^s |w|^{\alpha_n} w(t)  t^{N_{0,n}} dt\right) ds,  $$
          where $N_{0,n} = (N-1)(1+\alpha_n)$ 
           and   there exists some $\delta>0$ which does not depend on $n$,  such that  on $]r^\prime_n, r_n^\prime + \delta[$,    $\varphi_n^\prime >0$,  and $\varphi_n^{\prime\prime} >0$.

      Using remark \ref{remaru} there exists exactly one point $r^\prime$ such  that $\varphi$ is  decreasing  on  $]\bar r , r^\prime[$ and increasing on $]r^\prime,  \bar t[$, hence since $\varphi_n$ converges uniformly to $\varphi$,  one gets that $r_n^\prime$ converges to $r^\prime$.  
     
     From  the  definition of $T_n$ one sees that $\varphi_n$ converges uniformly  on $]r^\prime-{3\delta\over 4} , r^\prime[$ to the solution  $\psi$   on that intervall of 
         $T(\psi) = \psi$, 
         where 
         $$T(w)(r) = \varphi ( r^\prime )-\int_{r^\prime} ^r \varphi_{p^\prime} \left({\mu  (1+\alpha)\over A s^{N^-} }\int_{ r^\prime }^s |w|^{\alpha} w(t) t^{N^-} dt\right)ds.  $$
         
         This implies that $ \varphi$  is a ${\cal C}^2$ solution on $] r^\prime-{3\delta\over 4}, r^\prime[$ . 
          We do the same on  $]r^\prime, r^\prime +{ 3\delta\over 4} [$.

           We now consider  the equation on $]\bar r, r^\prime -{\delta\over 2}[$. As soon as $n$ is large enough in  order that $\bar r > r_{n-1}^\prime $,  on that intervall  $\varphi_n$ satisfies 
           $$(\varphi_n^\prime, \varphi_n^{\prime\prime}) = f_n (\varphi_n, \varphi_n^\prime)$$ 
           where 
           $f_n= (f_{1,n}, f_{2,n})$,
           $f_{1,n} (r, y_1, y_2) = y_2$, and  
              $$f_{2,n}(r, y_1, y_2) = M_n\left(-{m_n(y_2) (N-1)\over r }- {|y_1|^{\alpha_n} y_1 \over |y_2|^{\alpha_n}}\right), $$
       where $M_n$ and $m_n$ are respectively the functions 
       $$ M_n (x)= \left\{\begin{array}{c}
        {x\over A}\ {\rm if} \ x>0\\
         {x\over a_n} \ {\rm if } \ x<0,
         \end{array}\right.   $$
          and 
           $$ m_n (x)= \left\{\begin{array}{c}
        { A x }\ {\rm if} \ x>0\\
         {a_n x } \ {\rm if }\  x<0.
         \end{array}\right.$$

            It is clear that $f_n$ is uniformly Lipschitzian on $ ](\varphi(\bar r), \varphi (r^\prime)[\times ]\varphi^\prime (\bar r), \varphi^\prime(r^\prime -{\delta\over 2})[$. 
              Then $\varphi_n$ converges in ${\cal C}^1$ (even ${\cal C}^2$)  to some solution  $\psi$
              on $]\bar r, r^\prime -{\delta\over 2}[$ 
              of 
              $$(\psi^\prime, \psi^") = f(\psi, \psi^\prime)$$
              whith  
                 $f= (f_{1}, f_{2})$,
           $f_{1} (r, y_1, y_2) = y_2$,  
              $$f_{2}(r, y_1, y_2) = M\left(-{m(y_2) (N-1)\over r }- {|y_1|^{\alpha} y_1 \over |y_2|^{\alpha}}\right), $$
       and  $M$ and $m$ are respectively the functions 
       $$ M (x)= \left\{\begin{array}{c}
        {x\over A}\ {\rm if} \ x>0\\
         {x\over a} \ {\rm if } \ x<0, 
         \end{array}\right.$$
          and 
           $$ m(x)= \left\{\begin{array}{c}
        { A x }\ {\rm if} \ x>0\\
         {a x } \ {\rm if }\  x<0. 
         \end{array}\right.$$
          with the condition $\psi (\bar r)=0$, $\psi^\prime (\bar r) = \varphi^\prime (\bar r)$.
               
               This implies that $\varphi$ is  ${\cal C}^2$ on  $]\bar r, r^\prime -{\delta\over 2}[\cup [r^\prime -{3\delta\over 4}[ $. We can do the same on $]r^\prime + {\delta\over 2},  \bar t[$ and get in that way the regularity of $\varphi$ on $[r^\prime, \bar t[$.   In fact the proof  contains the regularity of $\varphi$ on a  open neighborhood of $[\bar r, \bar t]$. Since this can be repeated on each  intervall  delimited by two zeros of $\varphi$  one gets the regularity of $\varphi$ on $B(0,1)$.   As a consequence  of proposition \ref{propri} we have obtained that $\mu = \mu_k^+$. Since $\mu_k^+(\alpha_n, a_n)$ has a unique                         
cluster point we get that all the sequence converges to $\mu_k^+$.                       
           
    \section{Conclusion  and  supplementary results}

             Let $K_{\alpha,a } $ be the operator defined on ${\cal C}(\Omega)$  by : 
For  $f\in {\cal C} (\overline{\Omega})$, $K_{\alpha,a}(f)$ is  the unique $v\in {\cal C} (\overline{\Omega})$ solution of 
  $$\left\{\begin{array}{cc}
\tilde F_{\alpha,a}  (r,  v^\prime , v^") -  |v|^\alpha v =- f&  {\rm in}\  \Omega\\
  v = 0& \ {\rm on } \ \partial \Omega.
  \end{array}\right.
  $$
 The operator $K_{\alpha,a}$ is well defined since $\alpha>-1$, and defining   for $\mu$ positive  given $K_{\alpha, a, \mu} (u)=K_{\alpha,a} ( (\mu+1)|u|^\alpha u) $,  one can note that the fixed points of $K_{\alpha, a,  \mu} $ exist if  $\mu$ is an eigenvalue,  as some  eigenfunction associated.

We  will be able to derive from the continuity results in the last section  some results about the degree of the operator  $K_{\alpha, a, \mu} $ in function of the position of $\mu$ with respect to the eigengalues $\mu_k^\pm$.               
             Next we shall establish some bifurcation results for the equations  defined as follows 
               
               Let $f$ be defined as
 $(\mu, s)\mapsto f(\mu,s)$ which is "super-linear" in $s$ uniformly  with respect to $\mu$ in the sense that 
$$\lim_{s\rightarrow 0}\frac{f(\mu,s)}{|s|^{1+\alpha}}=0. $$
We also assume that $f$  is locally bounded and continuous in all its variables. 
 
Then we shall consider the problem
   \begin{equation}\label{prob}
   \left\{\begin{array}{cc}
   \tilde  F_{\alpha,a } (r,
    u^\prime, u^")   +  \mu |u|^\alpha u + f(\mu,u)=0&\ {\rm in} \ \Omega \\
    u=0 &{\rm on } \ \partial \Omega. 
    \end{array}\right.
 \end{equation} 
 for which we shall prove bifurcation results, completing the results already obtained in \cite{BDb}. 
 
 This will be the object of a forthcomming paper.


\begin{thebibliography}{99}
\bibitem{A} A. Anane, {\it Simplicit\'e et isolation de la premi\`ere
valeur propre du $p$-laplacien avec poids.} (French) [Simplicity and
isolation of the  eigenvalue of the $p$-Laplacian with weight]  C.
R. Acad. Sci. Paris SŽr. I Math.  305  (1987),  no. 16, 725--728



 \bibitem{BKJ} Belloni, M, B.  Kawohl, B.; Juutinen, P. {\it The $p$-Laplace eigenvalue problem as $p\to\infty$ in a Finsler metric.} J. Eur. Math. Soc. (JEMS) 8 (2006), no. 1, 123--138.

\bibitem{BNV} H. Berestycki, L.  Nirenberg, S.R.S. Varadhan, {\it  The principal 
eigenvalue and maximum principle for second-order elliptic operators in general domains.}
  Comm. Pure Appl. Math.  47  (1994),  no. 1, 47--92.







\bibitem{BD1}{ I. Birindelli, F. Demengel,} {\it Comparison principle and
Liouville type results for singular fully  nonlinear operators}, Ann. Fac. Sci Toulouse Math, (6)13 (2004), N.2, 261-287. .
\bibitem{BD2}{ I. Birindelli, F. Demengel,} {\it  Eigenvalue and
Maximum principle for fully nonlinear singular operators} Advances in Partial Diff. Equations.{\bf 11},1 (2006) 91-119.



\bibitem{BD3}{ I. Birindelli, F. Demengel,} 
{\it Eigenvalue, maximum principle and regularity for fully non linear
homogeneous operators} Comm. Pure and Applied Analysis {\bf 6} (2007).



\bibitem{BD4}{ I. Birindelli, F. Demengel,} 
{\it The Dirichlet problem for singular fully  nonlinear operators} Discrete Cont Dynamical systems (2007) special number, (2007), 110-121.
\bibitem{BDr}{ I. Birindelli, F. Demengel,} {\it Uniqueness of the first  eigenfunction for fully nonlinear equations: the radial case.},  To appear in ZAA, journal for mathematical Analysis

 
\bibitem{BD5} {I. Birindelli, F. Demengel},   {\it Eigenvalue and Dirichlet problem for fully-nonlinear operators in non smooth domains}, to appear in  Journal of Mathematical Analysis and its applications. 


\bibitem{BDb} {I. Birindelli, F. Demengel} {\it Bifurcation  for singular fully nonlinear equations}
    To appear in  "On the notions of solution to nonlinear elliptic problems: results and
 developments", Quaderni di Matematica, 23. Department of Mathematics,
 Seconda Universita  di Napoli, Caserta,  2008.

 
 
\bibitem{BDW}{ I. Birindelli, F. Demengel, J. Wigniolle,} {\it Strict
maximum principle}, Proceedings of Workshop on 
Second Order Subelliptic Equations and Applications
Cortona, (2003).

\bibitem{BEQ}{J. Busca,M.J. Esteban, A. Quaas} {\it Nonlinear eigenvalues and bifurcation problems for Pucci's operator} Annales de l'Institut H. Poincar\'e, Analyse
non-lin\'eaire, 22, (2005), no.  2, 187-206. 

\bibitem{CC} { L. Caffarelli, X. Cabr\'e,} Fully-nonlinear equations
Colloquium Publications 43, American Mathematical Society, Providence, RI,1995.


\bibitem{CGG} { Y. G. Chen, Y. Giga, S. Goto, } {\it Uniqueness and existence of viscosity solutions of generalized mean curvature flow equations.}  J. Differential Geom.  33  (1991),  no. 3, 749--786.
\bibitem{CIL} M.G. Crandall, H. Ishii, P.L. Lions, {\it User's guide to viscosity solutions of second order partial differential equations}  Bull. Amer. Math. Soc. (N.S.)  27  (1992),  no. 1, 1--67. 
\bibitem{C} Cuesta, Mabel {\it Eigenvalue problems for the $p$-Laplacian with indefinite weights. }
Electron. J. Differential Equations 2001, No. 33, 9 pp.

\bibitem{CT} Cuesta, Mabel; Tak\`a\v c, Peter {\it A strong comparison principle for positive solutions of degenerate elliptic equations.} Differential Integral Equations 13 (2000), no. 4-6, 721--746.
\bibitem{CL}{ A. Cutr\`i, F. Leoni,} {\it On the Liouville property for
fully-nonlinear equations} Annales de l'Institut H. Poincar\'e, Analyse
non-lin\'eaire, (2000), 219-245.

\bibitem{DM} Del Pino, Manasevich, {\it Global bifurcation from the eigenvalues of the $p$-Laplacian } 
Journal of Differential  equations, 92, n¡2, (1991), pp. 226-251. 
\bibitem{I}{H. Ishii} {\em Viscosity solutions of non-linear partial
differential equations } Sugaku Expositions vol 9 , (1996).


\bibitem{IL}{H. Ishii, P.L. Lions}, {\it Viscosity solutions of Fully- Nonlinear Second  Order Elliptic Partial Differential Equations} {  J. Differential Equations}  83  (1990),  no. 1, 26--78.

\bibitem{IY}{H. Ishii, Y. Yoshimura}, {\it Demi-eigen values for uniformly elliptic Isaacs op erators }, preprint. 

\bibitem{J}{P. Juutinen} {\it On the principal eigenvalue of a very badly degenerate equation},   J. Differential Equations (2007), 532--550.


\bibitem{PLL}  P.-L. Lions, \emph{ Bifurcation and optimal stochastic control}, Nonlinear Anal. \textbf{7} (1983), no. 2, 177--207. 

\bibitem{P} S. Patrizi {\it The Neumann problem for singular fully nonlinear operators}  J. Math. Pures et Appl., (9)  90, (2008), no. 3, 286-311.



\bibitem{LP2} Lucia, Marcello; Prashanth, S. {\it Simplicity of principal eigenvalue for $p$-Laplace operator with singular indefinite weight.} Arch. Math. (Basel) 86 (2006), no. 1, 79--89. 

\bibitem{OT} M. \^Otani, T. Teshima {\it On the first eigenvalue of some quasi linear elliptic equation} Proc. Japan Acad. Ser. A. Math. Sci. 64 (1988), no.1, 8--10.

\bibitem{Q} A. Quaas, {\it Existence of positive solutions to a ``semilinear'' equation involving the Pucci's operators in a convex domain}, submitted.

\bibitem{QS} A. Quaas, B. Sirakov, \emph{ On the principal eigenvalues and the Dirichlet problem for fully nonlinear operators.} C. R. Math. Acad. Sci. Paris \textbf{342} (2006), no. 2, 115--118.
\bibitem{QS} A. Quaas, B. Sirakov, \emph{ On the principal eigenvalues and the Dirichlet problem for fully nonlinear operators.}  Adv. math. 218, (2008), no. 1, 105-135. 

\end{thebibliography}
      \end{document}